\def\draft{0}

\documentclass[11pt]{article}
\usepackage[dvips,pagebackref,letterpaper=true,colorlinks=true,pdfpagemode=none,linkcolor=blue,citecolor=blue,pdfstartview=FitH]{hyperref}

\usepackage{comment}

\ifnum\draft=1
\newcommand{\mnote}[1]{[{\bf Madhur's Note: #1]}}
\newcommand{\lnote}[1]{[{\bf Luca's Note: #1]}}
\newcommand{\onote}[1]{[{\bf Omer's Note: #1]}}
\newcommand{\snote}[1]{[{\bf Salil's Note: #1]}}
\else
\newcommand{\mnote}[1]{}
\newcommand{\lnote}[1]{}
\newcommand{\onote}[1]{}
\newcommand{\snote}[1]{}
\fi

\usepackage{amsmath,amsfonts}

\newtheorem{theorem}{Theorem}[section]
\newtheorem{nt}{Selfnote}

\newtheorem{claim}[theorem]{Claim}

\newtheorem{remk}[theorem]{Remark}
\newtheorem{exmp}[theorem]{Example}

\def\FullBox{\hbox{\vrule width 8pt height 8pt depth 0pt}}

\def\qed{\ifmmode\qquad\FullBox\else{\unskip\nobreak\hfil
\penalty50\hskip1em\null\nobreak\hfil\FullBox
\parfillskip=0pt\finalhyphendemerits=0\endgraf}\fi}

\def\qedsketch{\ifmmode\Box\else{\unskip\nobreak\hfil
\penalty50\hskip1em\null\nobreak\hfil$\Box$
\parfillskip=0pt\finalhyphendemerits=0\endgraf}\fi}

\newenvironment{claimproof}{\begin{quotation} \noindent
{\bf Proof of claim:~~}}{\qedsketch\end{quotation}}

\newfont{\inhead}{eufm10 scaled\magstep1}

\newcommand{\calF}{{\cal F}}

\newcommand{\R}{{\mathbb R}}

\newcommand{\poly}{{\mathrm{poly}}}

\let\E=\Exp

\begin{document}

\title{New Proofs of the Green--Tao--Ziegler Dense Model Theorem:\\
An Exposition}

\author{
Omer Reingold\thanks{Faculty of Mathematics  and Computer Science, Weizmann Institute of Science, Rehovot 76100, Israel. {\tt
omer.reingold@weizmann.ac.il}. Research supported by US-Israel Binational Science Foundation grant 2006060.} \and
Luca Trevisan\thanks{Computer Science Division, U.C. Berkeley. \texttt{luca@cs.berkeley.edu}.
Work partly done while visiting
Princeton University and the IAS.  This material is based upon work supported by the National
Science Foundation under grants CCF-0515231 and CCF-0729137 and by the US-Israel Binational
Science Foundation under grant 2006060.} \and
Madhur Tulsiani\thanks{Computer Science Division, U.C. Berkeley. \texttt{madhurt@cs.berkeley.edu}
Work partly done while visiting
Princeton University. This material is based upon work supported by the National
Science Foundation under grants CCF-0515231 and CCF-0729137 and by the US-Israel Binational
Science Foundation under grant 2006060.} \and
Salil Vadhan\thanks{School of Engineering and Applied Sciences, Harvard University.
\texttt{salil@eecs.harvard.edu}.  Work done during a visit to U.C. Berkeley, supported by the Miller
Foundation for Basic Research in Science, a Guggenheim Fellowship,
US-Israel Binational Science Foundation grant 2006060,
and the Office of Naval Research grant N00014-04-1-0478.}}

\maketitle

\begin{abstract}
Green, Tao and Ziegler \cite{GT:primes,TZ06} prove ``Dense Model Theorems'' of the following form: if
$R$ is a (possibly very sparse) pseudorandom subset of set $X$, and $D$ is a dense subset of $R$,
then $D$ may be modeled by a set $M$ whose density inside $X$ is approximately the same
as the density of $D$ in $R$. More generally, they show that a function that is
majorized by a pseudorandom measure can be written as a sum of a bounded function
having the same expectation plus a function that is ``indistinguishable from zero.''
\snote{changed `more abstractly' to `more generally'.  Which version seems more abstract may depend on the
reader. rest of abstract substantially revised}  This theorem plays a key role in the proof of the Green--Tao Theorem~\cite{GT:primes} that
the primes contain arbitrarily long arithmetic progressions.

In this note, we present a new proof of the Green--Tao--Ziegler Dense Model Theorem, which was discovered
independently by ourselves~\cite{RTTV08} and Gowers~\cite{G08}.   Our presentation follows the argument in
\cite{RTTV08} (which in turn was inspired by Nisan's proof of the Impagliazzo Hardcore Set Theorem~\cite{I95}),
but is translated to the original notation of Green, Tao, and Ziegler.

We refer to our full paper \cite{RTTV08} for
variants of the result with connections and applications to computational complexity theory, and to Gowers'
paper \cite{G08} for applications of the proof technique to ``decomposition, ``structure,''
and ``transference''
theorems in arithmetic and extremal combinatorics (as well as a broader survey of such theorems).
\end{abstract}

\section{The Green-Tao-Ziegler Theorem}

Let $X$ be a finite universe. We use the notation $\E_{x\in X} f(x)
:= \frac {1}{|X|} \sum_{x\in X} f(x)$. For two functions $f,g: X\to \R$
we define their {\em inner product} as

\[ \langle f,g \rangle := \E_{x\in X} f(x)g(x) \]

A {\em measure} on $X$
is a function $g: X \to \R$ such that $g \geq 0$
and $\E_{x\in X} g(x) \leq 1$. A measure $g$
is {\em bounded} if $g \leq 1$.

Let $\cal F$ be a collection of bounded functions $f: X \to [-1,1]$.
We say that two measures $g,h$ are {\em $\epsilon$-indistinguishable}
according to $\cal F$ if

\[ \forall f\in \calF. | \langle g-h, f \rangle | \leq \epsilon \]

(It can be noted, although this fact will not be used, that
if we define $\|g\|_{\calF} = \max_{f\in \calF} | \langle g,f\rangle|$, then
$\|\cdot\|_{\calF}$ is a semi-norm, and we have that $g$ and $h$
are $\epsilon$-indistinguishable if and only if $\|g-h\|_{\calF} \leq \epsilon$.
Hence the notion of indistinguishability may be seen as a semi-metric
imposed on the space of functions $X \to \R$. If $\calF$ contains
{\em all} bounded functions $f:X\to [-1,1]$, then $\|\cdot\|_{\calF}$
is the standard $\ell_1$ norm.)

We say that a measure $g$ is {\em $\epsilon$-pseudorandom} according to $\cal F$
if $g$ and $1_X$ are {\em $\epsilon$-indistinguishable}
according to $\cal F$, where $1_X$ is the function that is identically equal to 1.

If $\cal F$ is a collection of bounded functions $f: X \to [-1,1]$,
we denote by ${\cal F}^k$ the collections of all functions of the form
$\prod_{i=1}^{k'} f_i$, where $f_i \in \cal F$ and $k' \leq k$.
\snote{changed $\calF_k$ to $\calF^k$}
In particular, if
$\cal F$ is closed under multiplication, then ${\cal F}^k = \cal F$.
\mnote{Changed definition of ${\cal F}^k$ slightly to have upto $k$ functions
instead of exactly $k$.}

\begin{theorem}[Green, Tao, Ziegler \cite{GT:primes,TZ06}]
For every $\epsilon>0$, there is a $k = (1/\epsilon)^{O(1)}$
and an $\epsilon' = \exp(-(1/\epsilon)^{O(1)})$ such that the following holds:

Suppose that $\cal F$ is a finite collection of bounded functions $f:X\to [-1,1]$ on a finite set $X$,
$\nu: X \to \R$ is an $\epsilon'$-pseudorandom measure
according to ${\cal F}^k$, and $g:X \to \R$ is a measure
such that $g \leq \nu$.

Then there is a bounded measure $g_1: X\to [0,1]$ such that
\begin{enumerate}
\item $\E_{x\in X} g_1(x) = \E_{x\in X} g(x)$, and \label{itm:density}
\item $g_1$ and $g$ are $\epsilon$-indistinguishable according to $\cal F$. \label{itm:indistinguishable}
\end{enumerate}
\end{theorem}
\snote{added finiteness condition on $X$; do we also need one on $\calF$, or does the min-max theorem
apply even when the space of strategies is infinite-dimensionsal?}

Green, Tao, and Ziegler~\cite{GT:primes,TZ06} state the conclusion in the following
equivalent form: we can write $g=g_1+g_2$,
where $g_1$ is a bounded measure, $g_1$ and $g$ have the same expectation, and $g_2$
is nearly orthogonal to $\cal F$ in the sense that
$| \langle g_2, f \rangle | \leq \epsilon$ for all $f\in\calF$.

\snote{following two pars are new}
We now describe how the theorem can be interpreted
as saying that ``every dense subset of a pseudorandom set has a dense model'',
as mentioned in the abstract.
From any sets $D\subseteq R\subseteq X$, we can obtain measures $\nu\geq g$ by setting
$\nu = 1_R \cdot |X|/|R|$ and $g = 1_D \cdot |X|/|R|$, where we write $1_S$ for the characteristic function of a set
$S$.  Then the condition that $\nu$ is $\epsilon'$-pseudorandom according to $\calF$
says that every function
$f\in \calF$ has the same average over $R$ as it does over $X$, to within $\pm \epsilon'$, which is a natural
pseudorandomness property of the set $R$.  And the expectation
of $g$ is precisely the density of $D$ in $R$, i.e. $|D|/|R|$.    Now, assuming that $R$ does indeed satisfy the
foregoing pseudorandomness property, let $g_1$ be the bounded function given in the conclusion of the theorem.
Suppose for
starters that $g_1$ is the characteristic function of some set $M\subseteq X$.  Then Item~\ref{itm:density} says
$M$ has the same density in $X$ as $D$ has in $R$.  And Item~\ref{itm:indistinguishable} says that $D$ and $M$ are
indistinguishable from each other, in the sense that every function in $\calF$ has the same average over both
sets, to within $\pm \epsilon/\delta$, where $\delta=|D|/|R|=|M|/|X|$.  So $M$ is indeed a ``dense model'' of $X$.

The actual theorem above can be interpreted as simply
allowing all of the sets, namely $D$ and $R$ in the hypothesis and $M$ in the conclusion, to have their
characteristic functions replaced with bounded measures of the same expectation.  We note that,
by an argument of Impagliazzo~\cite{I95}, allowing the
function $g_1$ in the conclusion to be a measure rather than the
characteristic function of some set $M$ does not substantially weaken the theorem.  Indeed, given $g_1$, we
can construct a set $M$ using the probabilistic method, including each element $x\in X$ in $M$ independently
with probability $g_1(x)$.  Then, by Chernoff Bounds,
$M$ will have density at least $(1-\epsilon)\delta$
and its characteristic function will be $2\epsilon$-indistinguishable from $g$
according to $\calF$ with probability $1-|\calF|\cdot \exp(-\Omega(\delta\epsilon^2|X|))$. \snote{calculations look right?}

\section{Our Proof}

We prove the contrapositive: assuming that $g_1$
is $\epsilon$-distinguishable from all dense models $g$ by functions in $\calF$, we prove
that $\nu$ cannot be pseudorandom, i.e. it is $\epsilon'$-distinguishable from $1_X$ by some
function in $\calF^k$.

Let $\delta:= \E_{x\in X} g(x)$ and let
us denote, for convenience, by $G$ the set of ``dense measures'' $g_1:X\to [0,1]$
such that $\E g_1 = \delta$. Our assumption can be written as

\[ \forall g_1 \in G. \exists f\in {\cal F}. | \langle g-g_1 , f \rangle | > \epsilon \]

If we denote by $\cal F'$ the closure of $\cal F$ under negation, that
is ${\cal F'} := {\cal F} \cup \{ -f : f\in {\cal F} \}$, we can remove the absolute values:

\begin{equation}
\label{eq:distinguish}
\forall g_1 \in G. \exists f\in {\cal F}'.  \langle g-g_1 , f \rangle  > \epsilon.
\end{equation}

\paragraph{Proof outline.}  Suppose that we can manage to find a $g_1,f$ pair such that the above holds
{\em and} for which $g_1(x)=1$ on every point in the support of $f$.
Then
it turns out that $f$ must also distinguish $\nu$ from $1_X$.  Indeed,
$\langle f,\nu\rangle \geq \langle f,g\rangle$, because $\nu\geq g$ pointwise,
and
$\langle f,1_X\rangle = \langle f,g_1\rangle$.

We will not be able to find such a $g_1,f$ pair with $f\in \calF'$, but we will be able to
do so with a function $f$ that is a {\em convex combination} of functions in $\calF'$ composed
with a {\em threshold function}.  Then we show how to convert a distinguisher of such a form into a
distinguisher that is a product of at most $k$ functions from $\calF$.

In more detail, the proof will proceed in the following steps:
\begin{enumerate}
\item By replacing $\calF'$ with its convex hull, we reverse the order of quantifiers
in (\ref{eq:distinguish}), and obtain a single $\bar{f}$ that $\epsilon$-distinguishes $g$ from {\em every} bounded measure $g_1$
of expectation $\delta$. \label{step:reverse}

\item With an appropriate choice of  $g_1$ (namely, the characteristic function of the $\delta|X|$ inputs on which
$\bar{f}$ is largest), we argue that a thresholded version of $\bar{f}$, denoted $\bar{f}_t$, continues to $\Omega(\epsilon)$-distinguish $g$ from $g_1$,
and has support contained in $g_1^{-1}(1)$.  By the above argument, $\bar{f}_t$ $\Omega(\epsilon)$-distinguishes
$\nu$ from $1_X$. \label{step:threshold}

\item By approximating the threshold function with a low-degree polynomial that has relatively small coefficients,
we deduce that there are at most $k$ functions from $\calF$ whose product $\epsilon'$-distinguishes $\nu$
from $1_X$. \label{step:polynomial}
\end{enumerate}

\paragraph{Proof Details.}
We now proceed with Step~\ref{step:reverse}, where we reverse the order of quantifiers in (\ref{eq:distinguish}).

\begin{claim}
There is a function $\overline f$ that is a convex combination of functions
from $\cal F'$ and satisfies:

\[ \forall g_1 \in G. \langle g-g_1 , \overline f \rangle  > \epsilon \]
\end{claim}

\begin{claimproof}
We use the min-max theorem for 2-player zero-sum games (which is a consequence of
the Hahn-Banach Theorem, as used in Gowers' version of the proof~\cite{G08}). \snote{confirm}
We think of a zero-sum game where the first player picks a function $f\in {{\cal F}'}$, the second player
picks a function $g_1 \in G$, and the payoff
is $\langle g-g_1 , f \rangle$ for the first player, and $-\langle g-g_1 , f \rangle$ for
the second player.

By the min-max theorem, the game has a ``value'' $\alpha$ for which the first player has
an optimal mixed strategy (a convex combination of strategies) $\bar f$,
and the second player has an optimal mixed strategy $\bar  g_1$,
such that
\begin{equation}\label{eq:minmax}
 \forall g_1 \in G, \ \  \langle g-g_1 , \bar f \rangle \geq \alpha \end{equation}
and
\begin{equation}\label{eq:minmaxb} \forall f \in {{\cal F}'}, \ \ \ \langle g-\bar g_1 , f \rangle \leq \alpha
\end{equation}
\mnote{Changed $\E \langle g-\bar g_1 , f \rangle$ to $\langle g-\bar g_1 , f \rangle$ in the above equation. It seems either both the equations should have the expectation or neither should, depending on how you interpret $\bar g_1$.}

Since $G$ is convex, $\bar g_1 \in G$, and our hypothesis
tells us that there exists a function $f$ such that
\[ \langle g-\bar g_1 , f \rangle > \epsilon \]

Taking this $f$ in
Inequality~(\ref{eq:minmaxb}), we get that $\alpha \geq \epsilon$.  The claim now follows from
Equation~(\ref{eq:minmax}).
\end{claimproof}

We now proceed with Step~\ref{step:threshold} of the proof.  Let
$S\subseteq X$ be the set of $\delta |X|$ elements of $X$
that maximize $\bar f$, and let $g_1$ be the characteristic function of $S$.\footnote{In case
$\delta|X|$ is not an integer, we define $g_1$ to be 1 on the $\lfloor \delta|X|\rfloor$ inputs
that maximize $\bar f$, to be 0 on the $|X|-\lceil \delta|X|\rceil$ inputs that minimize
$\bar f$, and to be an appropriate fractional value on the remaining element in order to
make the expectation of $g_1$ equal to $\delta$. \snote{new footnote}}
 \snote{replaced $1_S$ with
$g_1$ throughout}
Then $g_1$ is a bounded measure
of expectation $\delta$, i.e. an element of $G$, so we have:

\[  \langle g-g_1 , \bar f \rangle \geq \epsilon \ . \]

or, equivalently,

\begin{equation}
\label{eq:barf}
\langle g,\bar f \rangle \geq \langle g_1 , \bar f \rangle + \epsilon
\end{equation}

Now, we argue that by applying a threshold function to $\bar f$, we can ensure
that $g_1 = 1$ at every point in the support, while preserving the fact that
we distinguish $g$ from $g_1$.
Specifically, for a threshold $t$, define $\bar f_t : X\to \{0,1\}$ to be the boolean
function such that $\bar f_t(x) =1 $ if and only if $\bar f(x) \geq t$.
We will show that for some value of $t$, $\bar f_t$ has the properties we desire.
Moreover, it will be important for the final step to argue that the threshold is ``robust''
in the sense that it does not matter what happens in a small interval around $t$, where the discontinuity
of the threshold function could cause problems.  (Gowers~\cite{G08} handles this issue differently, by
instead showing that there is a distinguisher of the form $\max\{0,\bar{f}(x)-t\}$, which has the
advantage of being continuous everywhere as a function $\bar f(x)$.)

\begin{claim}\label{cl:threshold-disting}
There is a threshold $t\in \left[-1+\epsilon/3, 1 \right]$ such that

\[ \langle g,\bar f_t \rangle \geq
\left\langle g_1 , \bar f_{t-\epsilon/3} \right\rangle + \frac \epsilon 3\]
\end{claim}

\begin{claimproof}
First, observe that
\[ \bar f(x) =  \int_{-1}^1 \bar f_t(x) dt  - 1.\]
From (\ref{eq:barf}) and the fact that
$\langle g, 1_X \rangle = \langle g_1, {\mathbf 1 } \rangle = \delta$, we
have
\[ \langle g, \bar f+1 \rangle \geq \langle g_1, \bar f +1 \rangle + \epsilon, \]
which is equivalent to
\begin{equation} \label{eq:disting-integral}
\int_{-1}^1 \langle g,  \bar f_t \rangle dt \geq \int_{-1}^1 \langle g_1, \bar f_t  \rangle dt + \epsilon
\end{equation}
Now if the claim were false, we would have
\begin{eqnarray*}
\int_{-1}^1 \langle g,  \bar f_t \rangle dt
&=& \int_{-1}^{-1+\epsilon/3} \langle g,  \bar f_t \rangle dt
+ \int_{-1+\epsilon/3}^1 \langle g,  \bar f_t \rangle dt\\
& < & \int_{-1}^{-1+\epsilon/3} \langle g,  {\bf 1} \rangle dt
+ \int_{-1+\epsilon/3}^1
\left( \langle g_1,  \bar f_{t-\epsilon/3} \rangle + \frac \epsilon 3\right) dt\\
& \leq& \frac \epsilon 3 \cdot \delta + \int_{-1+\epsilon/3}^1
\langle g_1,  \bar f_{t-\epsilon/3} \rangle dt + \left( 2-\frac \epsilon 3 \right) \cdot \frac \epsilon 3 \\
& < & \int_{-1}^1
\langle g_1,  \bar f_{t} \rangle dt + \epsilon,
\end{eqnarray*}
contradicting Equation(\ref{eq:disting-integral}).
\end{claimproof}
\mnote{Changed $\bar f_{t-\frac \epsilon 3}$ to $\bar f_{t}$ in the last step of the proof, since the limits
of the integral are now different.}



We now argue that $g_1$ is identically equal to 1 on the support of $\bar{f}_{t-\epsilon/3}$.  Recall
$g_1$ is the characteristic function of the set of $\delta|X|$ inputs maximizing $\bar{f}$.  So if
$g_1(x)<1$ for some $x$ in the support of $\bar{f}_{t-\epsilon/3}$, then $g_1(x)=0$ everywhere outside
the support of $\bar{f}_{t-\epsilon/3}$.  \snote{rephrased this argument so it also applies for the case
that $\delta|X|$ is not an integer, and $g_1$ is non-boolean on a single point.}  But then
\[ \langle g_1, \bar f_{t-\epsilon/3}\rangle = \langle g_1,1_X \rangle
= \delta = \langle g, 1_X \rangle \geq \langle g,\bar f_t \rangle \]
in contradiction to Claim~\ref{cl:threshold-disting}.

Putting everything together, we have

\begin{eqnarray*}
\langle \nu, \bar f_t \rangle
&\geq& \langle g,\bar{f}_t\rangle\\
&\geq& \langle g_1,\bar{f}_{t-\epsilon/3}\rangle+\epsilon/3\\
&\geq& \langle 1_X,\bar{f}_{t-\epsilon/3}\rangle+\epsilon/3.
\end{eqnarray*}

Finally, we proceed with Step~\ref{step:polynomial}, where we
find a distinguisher that is defined as a product of
functions from $\cal F$, rather than being a threshold function
applied to a convex combination of elements of $\cal F'$.
We do this by approximating
the threshold function by a polynomial, using the following special case of the
Weierstrass Approximation Theorem.

\begin{claim}
For every $\alpha,\beta\in [0,1]$, $t\in [\alpha,1]$, there exists a
polynomial $p$
of degree $\poly(1/\alpha,1/\beta)$ and
with coefficients bounded in absolute value by ${\rm exp} (\poly(1/\alpha,1/\beta))$
such that
\begin{enumerate}
\item For all $z\in [-1,1]$, we have $p(z)\in [0,1]$.
\item For all $z\in [-1,t-\alpha]$, we have $p(z)\in [0,\beta]$.
\item For all $z\in [t,1]$, we have $p(z)\in [1 -\beta,1]$.
\end{enumerate}
\end{claim}

We set $\alpha=\epsilon/3$ and $\beta=\epsilon/12$ in the claim to obtain a polynomial
$p(z)=\sum_{i=0}^d c_iz^i$ of degree $d=\poly(1/\epsilon)$ with coefficients satisfying
$|c_i|\leq \exp(\poly(1/\epsilon))$ and such that for every $x$ we have
\[ \bar f_{t}(x) - \frac \epsilon {12} \leq (p\circ\bar f)(x)) \leq \bar f_{t-\epsilon/3 }(x) +
\frac \epsilon {12}, \]
where $\circ$ denotes composition.
From the properties of the polynomial $p$, we get
\[ \langle \nu, p\circ\bar{f} \rangle \geq \langle \nu, \bar f_t \rangle - \frac \epsilon {12} \]
and
\[ \langle 1_X, p\circ\bar{f} \rangle \leq \langle 1_X, \bar f_t \rangle + \frac \epsilon {12}, \]
giving
\begin{equation}
\langle \nu, p\circ\bar{f} \geq \langle 1_X, p\circ\bar{f} \rangle + \frac \epsilon 6.
\end{equation}

If the polynomial $p\circ\bar f = \sum_i c_i \bar f^i$ has inner product at least $\epsilon/6$
with $\nu -1$, there must  exist a single term $c_k \bar f^k$
whose inner product with $\nu-1$ is at least $\epsilon/(6(d+1))$, which in turn
implies that $\bar f^k$ has inner product of absolute value at least
$\epsilon':=\epsilon/(4(d+1)|c_k|) = \exp(-\poly(1/\epsilon))$ with $\nu-1$:
\[ | \langle \nu-1, \bar f^k \rangle| \geq \epsilon' \]
Suppose that
$\langle \nu-1, \bar f^k \rangle \geq \epsilon'$.  (The reasoning is analogous in the case of a negative inner
product.)
Recall
that the function $\bar f$ is a convex combination of functions from $\cal F'$.
\snote{removed: $\bar f(x) = \sum_{f\in \cal F'} \lambda_f f(x)$, since this assumes that
$\calF$ is finite}
This means that we may think
of $\bar f(x)$ as being the expectation of a random variable $f(x)$, in
which the function $f$ is picked according to some probability measure on
$\calF'$.
Then the value $\bar f(x)^k$ is the expectation of the process where
we sample independently $k$ functions $f_1,\ldots,f_k$ as before,
and then compute $\prod_i f_i (x)$. By linearity of expectation, we can
write
\[ \epsilon' \leq \langle \nu-1, \bar f(x)^k \rangle = \E \left\langle \nu-1, \prod_i f_i (\cdot) \right\rangle, \]
where the expectation is over the choices of the functions $f_i$ as described above.
We can now conclude that there is point in the sample space where a random
variable takes values at least as large as its expectation, and so
there are functions $f_1,\ldots,f_k \in \cal F'$ such that
\[ \left\langle \nu-1, \prod_i f_i (\cdot) \right\rangle \geq \epsilon' \]
Finally, replacing $f_i$ with $-f_i$ as appropriate, we can have all the functions $f_i$ be in $\calF$
itself (rather than $\calF'$).

\section*{Acknowledgments}
We thank Terence Tao, Avi Wigderson, Noga Alon, Russell Impagliazzo,
Yishay Mansour, and Timothy Gowers for comments, suggestions and references.


\end{document}